\documentclass[11pt]{article}
\usepackage{graphicx}
\usepackage{amsthm, amsmath, amssymb, tikz, bm}
\usepackage{mathtools}
\usepackage{xifthen}

\usepackage[margin=1.2in]{geometry} 

\usepackage{xcolor}
\usepackage{lmodern}
\usepackage[T1]{fontenc}
\usepackage[shortlabels]{enumitem}
\usepackage{todonotes}
\usepackage[colorlinks=true,
linkcolor=blue,citecolor=blue,
urlcolor=blue]{hyperref}

\usetikzlibrary{calc,shapes, backgrounds}

\newtheorem{lemma}{Lemma}
\newtheorem{theorem}{Theorem}

\newtheorem{conjecture}{Conjecture}

\newtheorem{case}{Case}

\title{Poset Ramsey Numbers for Boolean Lattices}

\author{
Linyuan Lu
\thanks{University of South Carolina, Columbia, SC 29208,
({\tt lu@email.sc.edu}).} \and
Joshua C. Thompson \thanks{University of South Carolina, Columbia, SC 29208,
({\tt joshuact@email.sc.edu}).} 
}

\begin{document}

\maketitle

\begin{abstract}
A subposet $Q'$ of a poset $Q$ is a \textit{copy of a poset} $P$ if there is a bijection $f$ between elements of $P$ and $Q'$ such that $x \le y$ in $P$ iff $f(x) \le f(y)$ in $Q'$. For posets $P, P'$, let the \textit{poset Ramsey number} $R(P,P')$ be the smallest $N$ such that no matter how the elements of the Boolean lattice $Q_N$ are colored red and blue, there is a copy of $P$ with all red elements or a copy of $P'$ with all blue elements. Axenovich and Walzer introduced this concept in \textit{Order} (2017), where they proved $R(Q_2, Q_n) \le 2n + 2$ and $R(Q_n, Q_m) \le mn + n + m$, where $Q_n$ is the Boolean lattice of dimension $n$. They later proved $2n \le R(Q_n, Q_n) \le n^2 + 2n$. Walzer later proved $R(Q_n, Q_n) \le n^2 + 1$. We provide some improved bounds for $R(Q_n, Q_m)$ for various $n,m \in \mathbb{N}$. In particular, we prove that $R(Q_n, Q_n) \le n^2 - n + 2$, $R(Q_2, Q_n) \le \frac{5}{3}n + 2$, and $R(Q_3, Q_n) \le  \frac{37}{16}n  + \frac{39}{16}$. We also prove that $R(Q_2,Q_3) = 5$, and $R(Q_m, Q_n) \le (m - 2 + \frac{9m - 9}{(2m - 3)(m + 1)})n + m + 3$ \text{ for all } $n \ge m \ge 4$.
\end{abstract}

\section{Introduction}
Ramsey theory roughly says that any $2$-coloring of elements in a sufficiently large discrete system contains a monochromatic system of given size.  In the domain of complete graphs, the classical Ramsey theorem states that for any two graphs $G$ and $H$ there is a integer $N_0$ such that if the edges of a complete graph $K_N$ with $N\geq N_0$ are colored in two colors then there exists either a red copy of $G$ or a blue copy of $H$ in $K_N$. The least such number $N_0$ is called the Ramsey number $R(G,H)$. This theorem was proved by Ramsey \cite{ramsey} in 1930, but the problem of exactly determining these,  ``multicolor'' Ramsey numbers, and $k$-uniform hypergraph Ramsey numbers remains open and is the subject of continuing research. For examples, see \cite{triple, stepping, hypergraph, sparse, cooley, triangles}.

In this paper, we will consider the poset Ramsey number instead of the graph Ramsey number. Given two posets $(P,\leq)$ and $(Q,\leq')$, we say $(P,\leq)$ is a {\em subposet} of $(Q,\leq')$, if there is an injective mapping $f\colon P\to Q$ such that for any $x, y \in P$ we have
\begin{equation}\label{eq:orderperserving}
    x\leq y \mbox{ if and only if } f(x)\leq' f(y).
\end{equation}
The image $f(P)$ is called a {\em copy} of $P$ in $Q$. A \textit{Boolean lattice} of dimension $n$, denoted $Q_n$, is the power set of an $n$-element ground set $X$ equipped with the inclusion relation. The \textit{2-dimension} of a poset $P$, defined by Trotter \cite{trotter} and denoted by $\dim_2(P)$, is the smallest $n$ such that $Q_n$ contains a copy of $P$.

A poset $X$ has \textit{Ramsey property} if for any poset $P$ there is a poset $Z$ such that when one colors the copies of $X$ in $Z$ red or blue, there is a copy of $P$ in $Z$ such that all copies of $X$ in this copy of $P$ are red or all of them are blue. The general problem of determining which posets have Ramsey property was solved by Ne$\breve{s}$et$\breve{r}$il and R$\ddot{o}$dl \cite{rodl}. In this paper, $X$ is the single-element poset. In other words, the elements of posets are colored instead of more complicated substructures.

For posets $P$ and $P'$, let the \textbf{poset Ramsey number} $R(P,P')$ be the least integer $N$ such that whenever the elements of $Q_N$ are colored in red or blue, there exists either a red copy of $P$ or a blue copy of $P'$. The focus of this paper is the case where $P$ and $P'$ are Boolean lattices $Q_m$ and $Q_n$ for $m,n \in \mathbb{N}$.
Axenovich and Walzer \cite{axenovich} give upper bound and lower bounds for $R(Q_m, Q_m)$ for various values of $m,n \in \mathbb{N}$. In particular, they prove the following.

\begin{theorem}\label{thm:axenovich}
For any integers $n,m \ge 1$,
\medskip

(i) $2n \le R(Q_n, Q_n) \le n^2 + 2n$,
\medskip

(ii) $R(Q_1, Q_n) = n + 1$,
\medskip

(iii) $R(Q_2, Q_n) \le 2n + 2$,
\medskip

(iv) $n + m \le R(Q_n, Q_m) \le mn + n + m$,
\medskip

(v) $R(Q_2, Q_2) = 4, R(Q_3, Q_3) \in \{ 7, 8 \}$,
\medskip

(vi) A Boolean lattice $Q_{3n \log (n)}$ whose elements are colored red or blue randomly and independently with equal probability contains a monochromatic copy of $Q_n$ asymptotically almost surely.
\end{theorem}

Walzer, in his master's thesis \cite{walzer}, improved the upper bound in Theorem \ref{thm:axenovich}, part (i) to the following.

\begin{theorem}\label{thm:walzer}
$R(Q_n, Q_n) \le n^2 + 1$.
\end{theorem}

Axenovich and Walzer also studied Ramsey numbers for Boolean algebras in \cite{axenovich}. A Boolean algebra $\mathcal{B}_n$ of dimension $n$ is a set system $\{ X_0 \cup \bigcup_{i \in I} X_i: I \subseteq [n] \}$, where $X_0, X_1, \dots, X_n$ are pairwise disjoint sets, $X_i \neq \emptyset$ for $i = 1, \dots, n$. Boolean algebras have a more restrictive structure than Boolean lattices. If a subset of $Q_N$ contains a Boolean algebra of dimension $n$, then it contains a copy of $Q_n$. The converse, however, is not always true. Gunderson, R$\ddot{o}$dl, and Sidorenko \cite{sidorenko} first considered the number $R_{Alg}(n)$, defined to be the smallest $N$ such that any red/blue coloring of subsets of $[N]$ contains a red or a blue Boolean algebra of dimesnion $n$. Here, "contains" means subset containment in $2^{[N]}$, not containment as a subposet. The following theorem states the best known bounds on $R_{Alg}(n)$. The lower bound is given without proof by Brown, Erd$\H{o}$s, Chung, and Graham \cite{brown}, and the upper bound was proved by Axenovich and Walzer \cite{axenovich}.

\begin{theorem}\label{algebra}
There is a positive constant $c$ such that
\[
2^{cn} \le R_{Alg}(n) \le \min\{ 2^{2^{n + 1} n \log n}, n R_h(K^n(2,\dots,2))\}.
\]
\end{theorem}

Here, $K^n(s,\dots,s)$ is a complete $n$-uniform $n$-partite hypergraph with parts of size $s$ and $R_h(K^n(2,\dots,2))$ is the smallest $N'$ such that any 2-coloring of $K^n(N',\dots,N')$ contains a monochromatic $K^n(2,\dots,2)$.

Gunderson, R$\ddot{o}$dl, and Sidorenko \cite{sidorenko} also considered the number $b(n,d)$, defined to be the maximum cardinality of a $\mathcal{B}_d$-free family contained in $2^{[n]}$. They proved the following bounds:

\[
n^{- \frac{(1 + o(1))d}{2^{d + 1} - 2}} \cdot 2^n \le b(n,d) \le 10^d 2^{-2^{1 - d}}d^{d - 
2^{-d}}n^{-1/2^d} \cdot 2^n.
\]

Johnston, Lu, and Milans \cite{milans} later used the Lubell function to improve the upper bound to the following, where $C$ is a constant:

\[
b(n,d) \le C n^{-1/2^d} \cdot 2^n.
\]
\smallskip

In this paper, we improve the upper bounds on the poset Ramsery numbers $R(Q_m, Q_n)$ given by Axenovich and Walzer in \cite{axenovich}. In Section \ref{mainresults}, we prove that for any integer $n \ge 1$,

\begin{theorem}\label{thm:2 and n}
$R(Q_2, Q_n) \le \frac{5}{3}n +2$.
\end{theorem}

\begin{theorem}\label{thm:n and n}
$R(Q_n, Q_n) \le n^2 - n + 2$.
\end{theorem}

\begin{theorem}\label{thm:3andn}
$R(Q_3, Q_n) \le  \frac{37}{16}n + \frac{39}{16}$.
\end{theorem}

\smallskip

In Section \ref{mainresults}, for all integers $n \ge m \ge 4$, we also prove the following.

\begin{theorem}\label{thm:m and n}
$R(Q_m, Q_n) \le (m - 2 + \frac{9m - 9}{(2m - 3)(m + 1)})n + m + 3$ \text{ for all } $n\geq m \ge 4$.
\end{theorem}

Additionally, we are now able to identify the following previously unknown poset Ramsey number.

\begin{theorem}\label{thm:2 and 3} $R(Q_2, Q_3) = 5$.
\end{theorem}

In Section \ref{notation}, we give more definitions and introduce notation. Also in Section \ref{notation}, we state and prove Lemma \ref{mapping}, the key embedding lemma we use to prove Theorems \ref{thm:2 and n}, \ref{thm:n and n}, \ref{thm:3andn}, and \ref{thm:m and n}. We prove theorems \ref{thm:2 and n}, \ref{thm:n and n}, \ref{thm:3andn}, \ref{thm:m and n}, and \ref{thm:2 and 3} in Section \ref{mainresults}, and we devote Section \ref{remarks} to concluding remarks.



\section{Notation and Key Lemma}\label{notation}

A \textit{partially ordered set}, or \textit{poset}, consists of a set $S$ together with a partial order $\leq$, which is a binary relation on $S$ satisfying
\begin{description}
\item [Reflexive Property:] $x\leq x$, for any $x\in S$.
\item [Transitive Property:] If $x\leq y$ and $y\leq z$ then $x\leq z$ for any $x,y,z\in S$.
\item [Antisymmetric Property:] If $x\leq y$ and $y\leq x$ then $x=y$ for any $x,y\in S$.
\end{description}

Let $[n]$ denote the set $\{1,2,\ldots, n\}$ and $Q_n=(2^{[n]}, \subseteq)$ be the poset over the family of all subsets  of $[n]$. The $k$-th level of $Q_n$ is the set of all $k$-element subsets of the ground set $[n]$, where $0 \le k \le n$.
For any two subsets (of $[n]$) $S\subset T$, let $Q_{[S,T]}$ be the induced poset of $Q_n$ over all sets $F$ such that $S\subseteq F \subseteq T$.
Let $Q^*_n:=Q_n\setminus \{\emptyset, [n]\}$.
Let $\hat{R}(Q_m, Q_n)$ denote the smallest $N$ such that any red/blue coloring of $Q^*_N$ contains either a red copy of $Q^*_m$ or a blue copy of $Q^*_n$. Equivalently, $\hat{R}(Q_m, Q_n)$ is the least $N$ such that if $\emptyset$ and $[N]$ are assumed to be both red and blue while the rest of $Q_N$ is colored either red or blue, then $Q_N$ contains either a red copy of $Q_m$ or a blue copy of $Q_n$.
For a subset $S \subseteq N$, let $\bar S$ denote the complement set of $S$ in $[N]$. When $S=\{x\}$, we simply write $\bar x$ for $\overline {\{x\}}$.
\bigskip

The following key lemma generalizes the blob lemma of Axenovich and Walzer (see \cite{axenovich}, Lemma 3).
The special case $a=b=0$ gives the blob lemma.

\begin{lemma}\label{mapping}
For any nonnegative integers $N$, $m$, $n$, $n'$, $a$, $b$ satisfying $N\geq n'\geq n\geq a+b$ and $N\geq m$, suppose that
the Boolean lattice $Q_N$ on the ground set $[N]$ is colored in two colors red and blue satisfying
\begin{enumerate}
    \item There is an injection $i\colon Q_n \to Q_{n'} \subset Q_N$ with the following properties.
    \begin{itemize}
    \item
    $i$ maps the bottom $a$-layers of $Q_n$ to blue sets.
    \item For all sets $S$ in the top $b$ layers of $Q_n$, $i(S) \cup ([N] \setminus [n'])$ is blue.
    \end{itemize}
    \item $N\geq n'+ (n+1-a-b)*m$.
\end{enumerate}
Then either a blue subposet $Q_n$ or a red subposet $Q_m$ exists in $Q_N$.
\end{lemma}

\noindent{\bf Proof of Lemma \ref{mapping}:}  Let $Q_N$ be the Boolean lattice on the ground set $[N]$ colored red and blue with the properties listed above.

Let $k = n + 1 - (a + b)$. Since $N \ge n' + (n + 1 - (a + b))*m = n+ k*m$, we can partition $[N]$ like so:

\[
[N] = [n'] \cup X_1 \cup X_2 \cup \dots \cup X_k
\]

where $|X_i| \ge m$ for all $i \in [k]$. With this partition in mind, we create an injection $f$ of $Q_n$ into the blue sets of $Q_N$. Consider the map $f: Q_n \to Q_N$ defined by

\[
f(\emptyset) = \emptyset
\]
\[
f(S) = i(S) \text{ for all } S \text{ with } |S| \le a
\]
\[
f(S) = i(S) \cup X_1^{\ast} \text{ for all } S \text{ with } |S| = a + 1
\]
\[
\vdots
\]
\[
f(S) = i(S) \cup X_1 \cup X_2 \cup \dots \cup X_{j}^{\ast} \text{ for all } S \text{ with } |S| = a + j
\]
\[
\vdots
\]
\[
f(S) = i(S) \cup X_1 \cup X_2 \cup \dots \cup X_{k}^{\ast} \text{ for all } S \text{ with } |S| = [n] - b
\]
\[
f(S) = i(S) \cup X_1 \cup X_2 \cup \dots X_k \text{ for all } S \text{ with } |S| \ge [n] - b + 1
\]
\[
f([n]) = [N].
\]

Here, $i(S) \cup X_1 \cup X_2 \cup \dots \cup X_{j}^{\ast}$ denotes an arbitrarily chosen blue element from the subposet with bottom element $S \cup X_1 \cup X_2 \cup \dots X_{j - 1} \cup \emptyset$ and top element $i(S) \cup X_1 \cup X_2 \cup \dots \cup X_{j - 1} \cup X_{j}$. If no such blue element exists, this entire subposet is red and $Q_N$ contains a red $Q_m$.

If such a blue element always exists,
this function is well-defined and preserves all the subset relations found in $Q_n$. Its image consists entirely of blue elements, so $Q_N$ contains a blue $Q_n$. \qed

\section{Proof of Theorems}\label{mainresults}

\begin{proof}[\textbf{Proof of Theorem \ref{thm:2 and n}}]
For any integer $n \geq 2$, let $N \in \mathbb{N}$ be such that there exists a red/blue coloring of $Q_N$ containing no red copy of $Q_2$ and no blue copy of $Q_n$. Consider such
a red-blue coloring $c$ of $Q_N$.  Let $T$ be a red element such that $\min \{ N - |T|, |T| \} \le \min \{ N - |T'| , |T'| \}$ for all red elements $T' \in Q_N$. Without loss of generality, let $N - |T| \le |T|$. Let $a := N - |T|$.
Let $S$ be a red element such that $|S| \le |S'|$ for all red elements $S' \in Q_{[\emptyset, T]}$. Let $b := |S|$.
\medskip

\noindent{\bf Claim a:} $|T| - |S| \le n + 1$. 

\noindent{\bf Proof of Claim a:} Proof by contradiction. Otherwise, suppose $|T| - |S| \geq n + 2$. Let $u,v$ be two red elements in $Q_{[S,T]}$. If $u$ and $v$ are incomparable, $\{ S, u, v, T \}$ form a red $Q_2$. So every red element in $Q_{[S,T]}$ lies on the same maximal chain. With the exception of this maximal chain, the rest of $Q_{[S,T]}$ is blue, and we can find a blue $Q_n$. \qed
\medskip

\noindent{\bf Claim b:} $N\leq 3n+1 -2(a+b)$.

\noindent{\bf Proof of Claim b:}
Otherwise, we assume $N\geq 3n+2-2(a+b)$.
We have $N \ge n + (n + 1 - (a + b))*2$. Since the bottom $a$-layers of $Q_N$ are all colored blue, the bottom $a$-layers of $Q_{[\emptyset,[n]]}$ are all colored in blue. If we let $m = 2$, by Lemma $\ref{mapping}$, $Q_N$ contains either a blue subposet $Q_n$ or a red subposet $Q_m$. \qed
\bigskip

From Claim a, we have
\begin{equation}\label{eq:2}
    a + b = N -(|T|-|S|) \ge N - (n+1).
\end{equation}
Combining \eqref{eq:2} with Claim b, we have
\begin{equation}
    N\leq 3n+1-2[N-(n+1)]=5n+3-2N.
\end{equation}
We get 
$$N\leq \frac{5n}{3}+1,$$

which gives us the desired result.

\end{proof}


\begin{proof}[\textbf{Proof of Theorem \ref{thm:n and n}}]
Let $n \in \mathbb{N}$. The result is known to hold for $n = 1$ and $n = 2$, so let $n \ge 3$. Let $\hat{R}(Q_n, Q_n)$ denote the smallest $N$ such that any red-blue coloring of $Q_N$, where $\emptyset$ and $[N]$ are assumed to be both red and blue, contains either a red or blue copy of $Q_n$. Equivalently, any red/blue coloring of $Q_N\setminus \{\emptyset, [N]\}$ contains either a red or blue copy of $Q_n^{\ast}$. To prove the theorem, we first prove the following claim.
\bigskip

\noindent\textbf{Claim c.}
$\hat{R}(Q_n, Q_n) \le n^2 - n$ for all $n \ge 3$.

\noindent\textbf{Proof of Claim c.} By way of contradiction, suppose there is a red-blue coloring $c$
of $Q_N$ (with $N=n^2-n$) such that $\emptyset$ and $[N]$ are colored both red and blue while all other elements of $Q_N$ only receive one color. Since $N = n^2 - n$, there are $n^2 - n \ge 2n$ singleton sets in the first row of $Q_N$. By the Pigeonhole Principle, there are $n$ sets in the first row of $Q_N$ with the same color. 
Without loss of generality, suppose at least $n$ of these sets are blue. Then there is a subposet $Q^{\ast}_n$ of $Q_N$ such that level 1 of $Q^{\ast}_n$ consists of some subset of these blue sets.

We consider an injection $i:Q_n \to Q^{\ast}_n \subset Q_N$, which maps the bottom $a = 2$ layers of $Q_n$ to blue sets. Also, we also consider the top $b = 1$ layer of $Q_N$ to be colored blue.
By Lemma \ref{mapping}, since $N \ge n^2 - n = n + (n - 2)*n = n + (n + 1 - a - b)*m$, either a blue subposet $Q_n$ or a red subposet $Q_m$ exists in $Q_N$. \qed

% ...





\bigskip

Let $N = n^2 - n + 2$. Consider a $Q_N$, and let $Q_N$ be colored with a coloring $c: Q_N \to \{ \text{ red, blue } \}$. We now consider the following cases.
\bigskip

\begin{case} Sets $\emptyset$ and $[N]$ are the same color. 
\end{case}

Without loss of generality, we assume both $\emptyset$ and $[N]$ are colored in red.
If we can find two blue sets $S$ and $T$ with $|S|=1$, $|T|=N-1$, and $S\subset T$,
then we can consider the $Q_{[S,T]}$. Since $|T|-|S|=N-2\geq n^2-n$, by Claim c,
$Q_{[S,T]}$ either a red or blue copy of $Q_n \setminus \{\emptyset, [n]\}$, which can be
extended to a red or blue copy of $Q_n$. 

If we fail to find such two blue sets $S$ and $T$, there are only three subcases:
\begin{enumerate}
    \item All level 1 sets are red.
    \item All level $N-1$ sets are red.
    \item There exists an element $x\in {N}$ such that $\{x\}$ an $[N]\setminus \{x\}$ are blue but all other sets in level 1 and level $N-1$ are red.
    \end{enumerate}

In subcase one, since $N\geq n+ n(n-2)$, we can partition $[N]=[n]\cup X_1\cup \cdots \cup X_{n-2}$ so that $|X_i|\geq n$. We map the first two layers and the last layer of $Q_n$ into $Q_N$ and extend this map as in the proof of Claim c to get a red copy of $Q_n$.
Subcase two is similar. In subcase three, similar argument works for the subposet
$Q_{[\emptyset, \bar x]}$, where $\bar x=[N]\setminus \{x\}$. Note that the first two layers of $Q_{[\emptyset, \bar x]}$ are red, while the top element $\bar x$ can be treated as red since $[N]$ is red.

\begin{case} Sets $\emptyset$ and $[N]$ are not the same color. Without loss of generality, suppose $\emptyset$ is red and $[N]$ is blue.
\end{case}

Suppose there is a pair $S,T$ of comparable elements, where $S$ is blue, $T$ is red, $|S| = 1$, and $|T| = N - 1$. Since $\emptyset$ is red and $S$ is blue, and $[N]$ is red and $T$ is blue, the poset $Q_{[S,T]}$ of dimension $n^2 - n$ can be viewed as having bottom and top elements colored both red and blue. By Claim c, $Q_{[S,T]}$ contains a red $Q_n$ or a blue $Q_n$.
\bigskip

Otherwise, there are only four subcases:
\begin{enumerate}
    \item All level 1 sets are red and all level $N-1$ sets are blue.
    \item All level 1 sets are red,  and there exists a red $N-1$-set.
    \item All level $N-1$ sets are blue, and there exists a blue $1$-set.
    \item There exists an element $x\in {N}$ such that all level 1 sets except $\{x\}$ are red and all level $N-1$ sets except $\bar x$ are blue.
    \end{enumerate}

A similar argument works for subcases 2, 3, and 4 since we can find a $Q_{[N-1]}$ so that there are three layers of one color. 

In subcase 1, suppose there exists a blue set in level $2$. Then we can find a blue $Q_{[N-2]}$ and a similar argument works. If there does not exist such a blue set, the bottom three layers of $Q_N$ are red.

In this case, since
\[n^2-n+2\geq n + (n-2)*n,\]
we can partition $[N] = [n] \cup X_1 \cup \cdots \cup X_{n - 2}$ so that $|X_i| \ge n$. We map the first three layers of $Q_n$ into $Q_N$ to get a red copy of $Q_n$. Applying Lemma \ref{mapping} with $a=3$ and $b=0$, we get the desired monochromatic copy of $Q_n$.





In any case where $N = n^2 - n + 2$, we have shown $Q_N$ must contain a red $Q_n$ or a blue $Q_n$. It follows that $R(Q_n, Q_n) \le n^2 - n + 2$, the desired result.

\end{proof}

\begin{proof}[\textbf{Proof of Theorem \ref{thm:3andn}}] For any integer $n \ge 4$, let $N \in \mathbb{N}$ be such that there exists a red/blue coloring of $Q_N$ containing no red copy of $Q_3$ and no blue copy of $Q_n$. Consider a red-blue coloring $c$ of $Q_N$. Let $T$ be a red element such that $\min \{ N - |T|, |T| \} \le \min \{ N - |T'| , |T'| \}$ for all red elements $T' \in Q_N$. Without loss of generality, let $N - |T| \le |T|$. Let $a := N - |T|$. Let $S$ be a red element such that $|S| \le |S'|$ for all red elements $S' \in Q_{ [ \emptyset, T ] }$. Let $b := |S|$. 


Let $\hat{R}(Q_3, Q_n)$ denote the smallest $N$ such that any red/blue coloring of $Q_N$, where $\emptyset$ and $[N]$ are assumed to be both red and blue, contains either a red copy of $Q_3$ or a blue copy of $Q_n$. Equivalently, any red-blue coloring of $Q_N^{\ast}$ contains either a red copy of $Q_3^{\ast}$ or a blue copy of $Q_n^{\ast}$. To prove the theorem, we first prove the following claim.
\medskip

\noindent\textbf{Claim d:} $\hat{R}(Q_3, Q_n) \le \frac{7}{4} n + \frac{9}{4}$ for all $n$.

\noindent{\textbf{Proof of Claim d:}} By way of contradiction, suppose there is a red-blue coloring $c$ of $Q^*_N$ (with $N \geq \frac{7}{4} n + \frac{9}{4}$) such that it contains neither red subposet $Q^*_3$ nor blue subposet $Q^*_n$.



Let $\ell = \lceil \frac{3}{8}n + \frac{5}{8} \rceil$ be a fixed integer. Consider the bottom $\ell$ layers of $Q_N$. We look for red sets $A_1, A_2, A_3$ with the following property.
\smallskip

\begin{equation}
\forall i \in [3], \exists x_i \in [N] \text{ such that } x_i \in A_i, \text{ but } x_i \not\in A_j \hspace{.1 in} \forall j \in [3] \backslash i. \hspace{.4 in} \label{xi}
\end{equation}
\medskip

We consider the following cases.
\medskip

\noindent\textbf{Case 1.} \textit{ There exist sets $A_1, A_2, A_3$ with property \ref{xi}.}
\bigskip


Since \[
\ell = \lceil \frac{3}{8}n + \frac{5}{8} \rceil = \lceil \frac{ \frac{7}{4}n - n + \frac{5}{4}}{2} \rceil
\le \lceil \frac{N - n - 1}{2} \rceil \le \frac{N - n + 1}{2} 
\]
\[
N + 1 \ge 2 \ell + n,
\]
\[
N + 1 \ge \ell + (\ell - 1) + n + 1
\]
\medskip

we are able to create an injection of $Q_3$ into the red sets of $Q_N$. Consider the map $f:Q_3 \to Q_N$ defined by
\medskip

\[
f(\{ i \}) = A_i \text{ for all } i \in [3],
\]
\[
f( \{ i, j \}) = X_{i,j}^{\ast} \text{ for all } \{i,j\} \subset [3].
\]
\medskip

Here, $X_{i,j}^{\ast}$ denotes an arbitrarily chosen red element from the subposet with bottom element $A_i \cup A_j$ and top element $\bar x_k$, where $\{ i, j, k \} = [3]$.
If no such red element exists, this entire $n$-dimenional subposet is blue and $Q_N$ contains a blue $Q_n$.

If such a red element always exists, this function is well-defined and preserves all the subset relations found in $Q_n$. Its image consists entirely of red elements, so $Q_N$ contains a red $Q_3$, a contradiction.

\bigskip

\noindent\textbf{Case 2.} \textit{There exist red sets $B_1, B_2, B_3$ in the top $\ell$ layers of $Q_N$ with the following property.}

\begin{equation}
\forall i \in [3], \exists x_i \in [N] \text{ such that } x_i \not\in B_i, \text{ but } x_i \in B_j \hspace{.1 in} \forall j \in [3] \backslash i. \hspace{.4 in} \label{nxi}
\end{equation}
\smallskip

This case is the same is as Case 1, except everything is flipped over the middle layer(s) of $Q_N$. Using a similar argument, we show that $Q_N$ contains a blue $Q_n$ or a red $Q_3$.







\bigskip

\noindent\textbf{Case 3.}
\textit{ There do not exist such sets $A_1, A_2, A_3$ or $B_1, B_2, B_3$}.
\bigskip

Suppose we are only able to find one red set $A_1$. Then every set of elements of $[N] \backslash A_1$ in the first $\ell$ layers is blue. Note that $|A_1| \le \ell - 1$. 

Suppose we are only able to find 2 sets with property $\ref{xi}$. Let $a_3$ be an arbitrarily chosen $\ell$-element subset of $A_1 \cup A_2$. We claim that every set of elements of $[N] \backslash a_3$ in the first $\ell$ layers is blue. Suppose this is not the case, and there is a red set $X \subseteq [N] \backslash a_3$ in the first $\ell$ layers. Since $| A_1 \cup A_2 | \le 2(\ell - 1)$, we know $| A_1 \cup A_2 \backslash a_3 | \le \ell - 2$. Thus, there exists an $x \in X$ such that $x \not\in A_1 \cup A_2$. We let $x$ be $x_3$, $X$ be $A_3$, and $A_1, A_2, A_3$ have property $\ref{xi}$, a contradiction. We can eliminate at most $\ell$ elements from $[N]$ and guarantee that sets formed from the remaining elements in the bottom $\ell$ layers are all blue.

Similarly, if we are only able to find at most two red sets with property $\ref{nxi}$, we can require the inclusion of at most $\ell$ elements from $[N]$ and guarantee that sets formed in the top $\ell$ layers of $Q_N$ are all blue.
Since $n < n + 1 = \frac{7}{4}n + \frac{9}{4} - 2(\frac{3}{8}n + \frac{5}{8})$
, we can define a mapping $i: Q_n \to Q^{\ast}_n \subset Q_N$ such that the bottom $\ell$ layers of $Q_n$ map to blue elements in the bottom $\ell$ layers of $Q_N$ and the top $\ell$ layers of $Q_n$ map to blue elements in the top $\ell$ layers of $Q_N$.
\smallskip




Since  \[ \ell = \lceil \frac{3}{8}n + \frac{5}{8} \rceil \]
\[
2 \ell = 2 \lceil \frac{3}{8}n + \frac{5}{8} \rceil \ge 2(\frac{3}{8}n + \frac{5}{8}) - 1
\]
\[
- 2 \ell \le - \frac{3}{4}n - \frac{1}{4}
\]
\[
1 - 2\ell = - \frac{3}{4}n + \frac{3}{4}
\]
\[
n + 1 - 2\ell = \frac{1}{4}n + \frac{3}{4}
= \frac{ \frac{3}{4}n + \frac{9}{4} }{3}
= \frac{ \frac{7}{4}n - n + \frac{9}{4}}{3} \le \frac{N - n}{3},
\]

we have

\[
N \ge n + (n + 1 - 2 \ell) * 3.
\]
\smallskip

The bottom $a = \ell$ layers and the top $b = \ell$ layers of $Q_n^{\ast}$ are blue and $m = 3$, by Lemma \ref{mapping}, $Q_N$ contains either a blue subposet $Q_n$ or a red subposet $Q_3$.
\medskip







In any case where $N \ge \frac{7}{4} n + \frac{9}{4}$ and $[N]$ and $\emptyset$ are colored both red and blue, we have shown that $Q_N$ must contain a red $Q_3$ or a blue $Q_n$. It follows that $\hat{R}(Q_3, Q_n) \le \frac{7}{4}n + \frac{9}{4}$. \qed
\bigskip

Suppose $a \neq 0$ and $b \neq 0$. It follows that $|T|-|S| + 1 \le \frac{7}{4} n + \frac{9}{4}$ for all $n \in \mathbb{N}$.
\medskip

\noindent\textbf{Claim e:} $ N \le n + 3 (n + 1 - (a + b)) - 1$.

\noindent{\textbf{Proof of Claim e:}} Otherwise, we assume $N \ge n + 3 (n + 1 - (a + b))$. Let $k = n + 1 - (a + b)$, so $N \ge a + b + 3k$.

We can partition $[N]$ like so:

\[
[N] = [n] \cup X_1 \cup X_2 \cup \dots \cup X_k,
\]

where $|X_i| \ge 3$ for all $i \in [k]$. With this partition in mind, we define a mapping $i: Q_n \to Q_n^{\ast} \subset Q_N$, an injection of $Q_n$ into the blue sets of $Q_N$. By Lemma \ref{mapping}, $Q_N$ contains either a blue copy of $Q_n$ or a red copy of $Q_3$, a contradiction. \qed


% ...




% ...


\bigskip

From Claim d, we have
\begin{equation}\label{eq:3}
    a + b = N -(|T|-|S|) \ge N - (\frac{7}{4}n + \frac{5}{4}).
\end{equation}
Combining \eqref{eq:3} with Claim e, we have
\begin{equation}
    N \le n + 3 (n + 1 - (a + b)) - 1 \le n + 3(n + 1 -(N - \frac{7}{4}n - \frac{5}{4})) - 1.
\end{equation}
We get 
$$N \leq \frac{37}{16}n + \frac{23}{16}.$$
\bigskip


% ...





Now suppose $a = 0$. We consider the remaining two cases. In each case, we assume, by way of contradiction, that $N > \frac{37}{16}n + \frac{23}{16}.$

\bigskip

\noindent\textbf{Case 1.} \textit{ $a = 0 \text{ and } b = 0$. }
\medskip

In this case, both $\emptyset$ and $[N]$ are necessarily red. 
If we can find two blue sets $S$ and $T$ with $|S| = 1$, $|T| = N - 1$, and $S \subset T$, then we can consider $Q_{[S,T]}$. In this case, since $\emptyset$ is red and $S$ is blue, we can consider the bottom element of $Q_{[S,T]}$ to be both red and blue. Since $[N]$ is red and $T$ is blue, we can consider the top element of $Q_{[S,T]}$ to be both red and blue. By Claim d, $\hat{R}(Q_3, Q_n) + 2 \le \frac{7}{4}n + \frac{9}{4} + 2 < \frac{37}{16}n + \frac{23}{16}$ for $n \ge 4$.

If we cannot find such sets $S$ and $T$, we are left with the following three subcases:
\begin{enumerate}
    \item All level 1 sets are red.
    \item All level $N-1$ sets are red.
    \item There exists an element $x\in {N}$ such that $\{x\}$ and $\bar x$ are blue but all other sets in levels and $N-1$ are red.
    \end{enumerate}
    
In subcase 1, since $N > \frac{37}{16}n + \frac{23}{16} > n + 3 = 3 + (3 + 1 - 1 - 2) \ast n$, we can partition $[N] = [3] \cup X$ with $|X| \ge n$. We map the first $b=2$ layers and the last $a=1$ layer of $Q_3$ into $Q_N$. By Lemma 1, $Q_N$ contains either a blue $Q_n$ or a red $Q_3$. Subcase 2 is similar.

In subcase 3, a similar argument works for $Q_{[\emptyset, \bar x]}$. Note that the first two layers of $Q_{[\emptyset, \bar x]}$ are red, while the top element $\bar x$ can be treated as red since $[N]$ is red.

\medskip

\noindent\textbf{Case 2.} \textit{ $a = 0 \text{ and } b \neq 0$}.
\medskip

In this case, $\emptyset$ is necessarily blue and $[N]$ is necessarily red. Suppose there is a pair $S,T$ of comparable elements, where $S$ is red, $T$ is blue, $|S| = 1$, and $|T| = N - 1$. Since $\emptyset$ is blue and $S$ is red, and $T$ is blue and $[N]$ is red, the poset $Q_{[S,T]}$ of dimension at least $N- 2> \frac{37}{16}n + \frac{23}{16} - 2 > \frac{7}{4}n$ can be viewed as having bottom and top elements colored both red and blue. By Claim d, $Q_{[S,T]}$ contains a red $Q_3$ or a blue $Q_n$.

\medskip

Otherwise, there are only four subcases:
\begin{enumerate}
    \item All level 1 sets are blue.
    \item All level $N-1$ sets are red, and there exists a red $1$-set.
    \item There exists an element $x\in {N}$ such that all level 1 sets except $\{x\}$ are blue and all level $N-1$ sets except $\bar x$ are red.
    \end{enumerate}
    
In subcase 2, let $S$ be the red 1-set. Since $N - 1 > \frac{37}{16}n + \frac{23}{16} - 1 > n + 3 = 3 + (3 + 1 - 1 - 2) \ast n$, we can map the first $b = 1$ layer and the last $a = 2$ layers of $Q_3$ into $Q_{[S,[N]]}$. By Lemma 1, $Q_N$ contains either a blue $Q_n$ or a red $Q_3$.

In subcase 3, a similar argument works for $Q_{[x, [N]]}$. Note that the bottom layer and top two layers of $Q_{[x, [N]]}$ are red.

In subcase 1, let $S$ be a red set such that $|S| \le |S'|$ for all red sets $S'$ in $Q_N$. Suppose $S$ is in level $\ell$. Suppose $\ell \ge n + 1$. Then the bottom $n + 1$ layers of $Q_N$ are blue, and $Q_N$ contains a blue copy of $Q_n$. Thus $\ell \le n$.

Suppose there exists a blue $N-1$- set $T$. Suppose $\ell \le \lfloor \frac{9}{16}n - \frac{1}{4} \rfloor$. Since $\emptyset$ is blue and $S$ is red, and $T$ is blue and $[N]$ is red, we consider the top and bottom elements of $Q_{[S,T]}$ to be both red and blue. We have $N + 1 > \frac{37}{16}n + \frac{23}{16} + 1 > \frac{9}{16}n - \frac{1}{4} + 1 + \frac{7}{4}n \ge \ell + 1 + \frac{7}{4}n$. By Claim d, $Q_{[S,T]}$ contains a red $Q_3$ or a blue $Q_n$, so $Q_N$ contains a red $Q_3$ or a blue $Q_n$.

Suppose there exists a blue $N-1$-set $T$, but $\ell \ge \lfloor \frac{9}{16}n + \frac{3}{4} \rfloor$. We have $$N > \frac{37}{16}n + \frac{23}{16} > n + (\frac{7}{16}n - \frac{1}{4})\ast 3 = n + (n + 1 - 1 - \frac{9}{16}n - \frac{1}{4}) \ast 3 $$
$$
\ge n + (n + 1 - 1 - \ell )\ast 3.
$$

We can map the first $b = \ell$ layers and the last $a = 1$ layer of $Q_n$ into $Q_N$. By Lemma 1, $Q_N$ contains either a blue $Q_n$ or a red $Q_3$.

Now suppose there is no blue $N-1$-set. That is, the top 2 layers of $Q_N$ are red. We consider $Q_{[S,[N]]}$. Since $N > \frac{37}{16}n + \frac{23}{16} \ge 2n + 3 = \ell + 3 + (3 + 1 - 1 - 2)\ast n$, we map the first $b = 1$ layer and the last $a=2$ layers of $Q_3$ into $Q_{[S,[N]]}$. By Lemma 1, $Q_{[S,[N]]}$ contains either a blue $Q_n$ or a red $Q_3$, so $Q_N$ contains either a blue $Q_n$ or a red $Q_3$.
    



\end{proof}

\begin{proof}[\textbf{Proof of Theorem \ref{thm:m and n}}] For any integers $m,n \in \mathbb{N}$ with $n \ge m \ge 3$, let $N \in \mathbb{N}$ be such that there exists a red/blue coloring of $Q_N$ containing no red copy of $Q_3$ and no blue copy of $Q_n$. Consider a red-blue coloring $c$ of $Q_N$. Let $T$ be a red element such that $\min \{ N - |T|, |T| \} \le \min \{ N - |T'| , |T'| \}$ for all red elements $T' \in Q_N$. Without loss of generality, let $N - |T| \le |T|$. Let $a := N - |T|$. Let $S$ be a red element such that $|S| \le |S'|$ for all red elements $S' \in Q_{ [ \emptyset, T ] }$. Let $b := |S|$.


Let $\hat{R}(Q_m, Q_n)$ denote the smallest $N$ such that any red/blue coloring of $Q_N$, where $\emptyset$ and $[N]$ are colored both red and blue, contains either a red copy of $Q_m$ or a blue copy of $Q_n$. Equivalently, any red-blue coloring of $Q_N^{\ast}$ contains either a red copy of $Q_m^{\ast}$ or a blue copy of $Q_n^{\ast}$. To prove the theorem, we first prove the following claim.
\bigskip

\noindent\textbf{Claim f:}
$\hat{R}(Q_m, Q_n) \le (m - 2 + \frac{3}{2m - 3})n + m$ for all $n \ge m \ge 4$.

\noindent\textbf{Proof of Claim f:} By way of contradiction, suppose there is a red-blue coloring $c$ of $Q_N^{\ast}$ (with $N = (m - 2 + \frac{3}{2m - 3})n + m)$ such that it contains neither red subposet $Q_m^{\ast}$ nor blue subposet $Q_n^{\ast}$.

Let $\ell = \lceil 1 + \frac{3n}{m (2m - 3)} \rceil$ be a fixed integer. Consider the bottom $\ell$ layers of $Q_N$. We look for red sets $A_1, A_2, \dots, A_m$ with the following property.



\begin{equation}
\forall i \in [m], \exists x_i \in [N] \text{ such that } x_i \in A_i, \text{ but } x_i \not\in A_j \hspace{.1 in} \forall j \in [m] \backslash i. \hspace{.4 in} \label{mxi}
\end{equation}
\medskip



We consider the following cases.
\bigskip

\noindent\textbf{Case 1.} \textit{ There exist sets $A_1, A_2, \dots, A_m$ with property \ref{mxi}.}
\bigskip

Since\[
\ell = \left\lceil 1 + \frac{3n}{m (2m - 3)} \right\rceil
\]
\[
\ell \le 2 + \frac{3n}{m (2m - 3)}
\]
\[
\ell - 1 \le 1 + \frac{3n}{m (2m - 3)}
\]
\[
m( \ell - 1) \le m + \frac{3n}{2m - 3}
\]
\[
m ( \ell - 1) + n(m - 2) + 1 \le m + \frac{3n}{2m - 3} + n(m - 2) + 1
\]
\[
m ( \ell - 1) + n(m - 2) + 1 \le N + 1,
\]
\smallskip

we are able to create an injection of $Q_m$ into the red sets of $Q_N$. We can partition $[N]$ like so:

\[
[N] = [n] \cup X_1 \cup X_2 \cup \dots \cup X_{m - 1},
\]

where $X_1 = \bigcup_{i = 1}^m (A_i \backslash \{ x_i \})$ and $|X_i| \ge n$ for all $i$ with $2 \le i \le m - 1$. We create an injection of $Q_m$ into the red sets of $Q_N$. Consider the map $f:Q_m \to Q_N$ defined by
\medskip

\[
f(\emptyset) = \emptyset
\]
\[
f(\{ i \}) = A_i \text{ for all } i \in [m]
\]
\[
f( \{ i, j \}) = A_i \cup A_j \cup X_2^{\ast} \text{ for all } \{i,j\} \subset [m]
\]
\[
\vdots
\]
\[
f(S) = \bigcup_{i \in S} A_i \cup X_2 \cup \dots \cup X_d^{\ast} \text{ for all } S \subset [m] \text{ with } |S| = d
\]
\[
\vdots
\]
\[
f([m]) = [N].
\]
\smallskip

Here, $\bigcup_{i \in S} A_i \cup X_2 \cup \dots \cup X_d^{\ast}$ denotes an arbitrarily chosen red element from the subposet with bottom element $\bigcup_{i \in S} A_i \cup X_2 \cup \dots \cup X_{d-1}$ and top element $\bigcup_{i \in S} A_i \cup X_2 \cup \dots \cup X_d$.
If no such red element exists, this entire $n$-dimenional subposet is blue and $Q_N$ contains a blue $Q_n$.

If such a red element always exists, this function is well-defined and preserves all the subset relations found in $Q_n$. Its image consists entirely of red elements, so $Q_N$ contains a red $Q_m$.

\bigskip

\noindent\textbf{Case 2.} \textit{There exist red sets $B_1, B_2, \dots, B_m$ in the top $\ell$ layers of $Q_N$ with the following property.}

\begin{equation}
\forall i \in [m], \exists x_i \in [N] \text{ such that } x_i \not\in B_i, \text{ but } x_i \in B_j \hspace{.1 in} \forall j \in [m] \backslash i. \hspace{.4 in} \label{mnxi}
\end{equation}
\medskip

This case is the same as Case 1, except everything is flipped over the middle layer(s) of $Q_N$. Using a similar argument, we show that $Q_N$ contains a blue $Q_n$ or a red $Q_3$.

\bigskip

\noindent\textbf{Case 3.}
\textit{ There do not exist such sets $A_1, A_2, \dots, A_m$ or $B_1, B_2, \dots, B_m$}.
\bigskip

Suppose we are only able to find at most most $m-1$ sets $A_1, A_2, \dots, A_{m-1}$ with property \ref{mxi}. Let $a_m$ be an arbitrarily chosen subset of $\bigcup_{i = 1}^{m-1} A_i$ such that $|a_m| = (m-2)(\ell - 1) + 1$. We claim that every set of elements of $[N] \backslash a_m$ in the first $\ell$ layers is blue. Suppose this is not the case, and there is a red set $X \subseteq [N] \backslash a_m$ in the first $\ell$ layers.
Since $|\bigcup_{i = 1}^{m-1} A_i| \le (m - 1)(\ell - 1)$, we know $|\bigcup_{i = 1}^{m-1} A_i \backslash a_m| \le \ell - 2$. Thus, there exists an $x \in X$ such that $x \not\in \bigcup_{i = 1}^{m-1} A_i$. We let $x$ be $x_m$, $X$ be $A_m$, and $A_1, A_2, \dots, A_m$ have property \ref{mxi}, a contradiction. We can eliminate at most $(m - 2)(\ell - 1) + 1$ elements from $[N]$ and guarantee that sets formed from the remaining elements in the bottom $\ell$ layers are all blue.

Similarly, if we are only able to find at most $m-1$ red sets with property \ref{mnxi}, we can require the inclusion of at most $(m - 2)(\ell - 1)$ elements from $[N]$ and guarantee that sets formed in the top $\ell$ layers of $Q_N$ are all blue.

Since $n < N - 2 (m - 2)(\ell - 1)$ for all $m,n \ge 4$, we can define a mapping $i: Q_n \to Q^{\ast}_n \subset Q_N$ such that the bottom $\ell$ layers of $Q_n$ map to blue elements in the bottom $\ell$ layers of $Q_N$ and the top $\ell$ layers of $Q_n$ map to blue elements in the top $\ell$ layers of $Q_N$.

Since  \[ \ell = \lceil 1 + \frac{3n}{m(2m - 3)} \rceil \]
\[
\ell \ge 1 + \frac{3n}{m(2m - 3)}
\]
\[
\ell - 1 \ge \frac{3n}{m(2m - 3)}
\]
\[
(m - 2)(\ell - 1) \ge \frac{3n (m-2)}{m(2m - 3)}
\]
\[
-2(m - 2)(\ell - 1) \le \frac{-6n (m-2)}{m(2m - 3)}
\]
\[
n + 1 -2(m - 2)(\ell - 1) \le n + 1 + \frac{-6n (m-2)}{m(2m - 3)}
\]
\[
= \frac{(m - 3 + \frac{3}{2m - 3})n + m}{m}
\]
\[
= \frac{(m - 2 + \frac{3}{2m - 3})n  + m - n}{m} \le \frac{N - n}{m},
\]

we have

\[
N \ge n + (n + 1 - 2(m - 2)(\ell - 1)) * m.
\]
\smallskip

The bottom $a = \ell$ layers and the top $b = \ell$ layers of $Q_n^{\ast}$ are blue, By Lemma \ref{mapping}, $Q_N$ contains either a blue subposet $Q_n$ or a red subposet $Q_m$.
\medskip

In any case where $N \ge (m - 2 + \frac{3}{2m - 3})n + m$ and $[N]$ and $\emptyset$ are colored both red and blue, we have shown that $Q_N$ must contain a red $Q_m$ or a blue $Q_n$. It follows that $\hat{R}(Q_3, Q_n) \le (m - 2 + \frac{3}{2m - 3})n + m$. \qed

\bigskip

Suppose $a \neq 0$ and $b \neq 0$. It follows that $|T| - |S| + 1 \le (m - 2 + \frac{3}{2m - 3})n + m$ for all $n \ge m \ge 4$.
\medskip

\noindent\textbf{Claim g:} $N \le n + m(n + 1 - (a + b)) - 1$

\noindent\textbf{Proof of Claim g:} Otherwise, we assume $N \ge n + m(n + 1 - (a + b))$. Let $k = n + 1 - (a + b)$, so $N + 1 \ge a + b + mk$.

We can partition $[N]$ like so:

\[
[N] = [n] \cup X_1 \cup X_2 \cup \dots \cup X_k,
\]

where $k = \frac{N - n}{m}$ and $|X_i| \ge m$ for all $i \in [k]$. With this partition in mind, we define a mapping $i:Q_n \to Q_n^{\ast} \subset Q_N$, an injection of $Q_n$ into the blue sets of $Q_N$. By Lemma 1, $Q_N$ contains either a blue copy of $Q_n$ or a red copy of $Q_m$, a contradiction. \qed

% ...






% ...

\bigskip

From Claim f, we have

\begin{equation}\label{eq:60}
    a + b = N -(|T|-|S|) \ge N - ((m - 2 + \frac{3}{2m-3})n + m - 1).
\end{equation}
Combining \eqref{eq:60} with Claim g, we have
\[
    N \le n + m( n + 1 - (a + b)) - 1
\]
\begin{equation}
    \le n + m(n + 1 - (N - (m - 2 + \frac{3}{2m - 3})n - m + 1)) - 1
\end{equation}
We get
\begin{equation}\label{eq:blueblue}
N\leq (m - 2 + \frac{9m - 9}{(2m - 3)(m + 1)})n + m -1.
\end{equation}
\bigskip

Now suppose $a = 0$. We consider the remaining two cases. In each case, we assume, by way of contradiction, that $N > (m - 2 + \frac{9m - 9 }{(2m - 3)(m + 1)})n + m + 2$.
\bigskip

\noindent\textbf{Case 1.} \textit{ $a = 0 \text{ and } b = 0$. }
\medskip

In this case, both $\emptyset$ and $[N]$ are necessarily red. We consider levels 1, 2, $N - 2$, and $N - 1$. 
If we can find two blue sets $S$ and $T$ with $|S| \le 2$, $|T| \ge N - 2$, $|T| - |S| \ge N - 3$, and $S \subset T$, then we can consider $Q_{[S,T]}$. In this case, since $\emptyset$ is red and $S$ is blue, we can consider the bottom element of $Q_{[S,T]}$ to be both red and blue. Since $[N]$ is red and $T$ is blue, we can consider the top element of $Q_{[S,T]}$ to be both red and blue. By Claim f, $\hat{R}(Q_m, Q_n) + 3 \le (m - 2 + \frac{3}{2m - 3})n + m + 3 < (m - 2 + \frac{9m - 9}{(2m - 3)(m + 1)})n + m + 2 < N$ for sufficiently large $m$ and $n$.

If we cannot find such sets $S$ and $T$, we are left with the following subcases:
\begin{enumerate}
    \item All sets in levels 1 and 2 are red.
    \item All sets in levels $N-2$ and $N- 1$ are red.
    \item All sets in levels 1 and $N - 1$ are red.
    \item There exist blue sets $S$ and $T$ with $|S| \le 2$ and $|T| \ge N - 2$, but $S \not\subset T$.
    \end{enumerate}
    
In subcase 1, since $N > (m - 2 + \frac{9m - 9 }{(2m - 3)(m + 1)})n + m + 2 > m + (m - 3)\ast n$, we can partition $[N] = [m] \cup X_1 \cup X_2 \cup \dots \cup X_{m - 3}$ with $|X_i| \ge n$ for all $i \in [m - 3]$. We map the first $b = 3$ layers and the last $a = 1$ layer of $Q_m$ into $Q_N$. By Lemma 1, $Q_N$ contains either a blue $Q_n$ or a red $Q_m$. Subcase 2 is similar. Subcase 3 is similar, except we map the first $b = 2$ layers and the last $a = 2$ layers of $Q_m$ into $Q_N$.

In subcase 4, either $S$ is in level 1 or $T$ is in level $N - 1$. Otherwise, we apply the same strategy as in subcases 1, 2, or 3. Suppose, without loss of generality, that $T$ is in level $N - 1$. Then a similar argument works for $Q_{[\emptyset, T]}$. Note that the first three layers of $Q_{[\emptyset, T]}$ are red, while the top element $T$ can be treated as red since $[N]$ is red.

\medskip

\noindent\textbf{Case 2.} \textit{ $a = 0 \text{ and } b \neq 0$}.
\medskip

In this case, $\emptyset$ is necessarily blue and $[N]$ is necessarily red. Suppose there is a pair $S,T$ of comparable elements, where $S$ is red, $T$ is blue, $|S| \le 2$, $|T| \ge N - 2$, and $|T| - |S| \ge N - 4$. Since $\emptyset$ is blue and $S$ is red, and $T$ is blue and $[N]$ is red, we can consider the top and bottom elements of $Q_{[S,T]}$ to be both red and blue. 
By Claim f, $\hat{R}(Q_m, Q_n) + 4 \le (m - 2 + \frac{3}{2m - 3})n + \frac{m}{2m - 3} + 4 < (m - 2 + \frac{9m - 9}{(2m - 3)(m + 1)})n + m + 2$ for sufficiently large $m$ and $n$.

\medskip

Otherwise, there are only four remaining subcases:
\begin{enumerate}
    \item All sets in levels 1 and 2 are blue and all sets in levels $N-2$ and $N-1$ are red.
    \item All sets in levels 1 and 2 are blue and there exists a blue set $T$ with $|T| \ge N-2$.
    \item All sets in levels $N - 2$ and $N - 1$ are red and there exists a red set $S$ with $|S| \le 2$.
    \item There exists a red set $S$ and a blue set $T$ with $|S| \le 2$ and $|T| \ge N - 2$, but $S \not\subset T$.
\end{enumerate}

A similar argument works for subcases 3 and 4 since we can find a $Q_{[N-2]}$ so that there are four red layers.
In subcase 3, we consider $Q_{[S,[N]]}$. Since $N - 2 > (m - 2 + \frac{9m - 9}{(2m - 3)(m + 1)})n + m + 2 - 2 > m + (m - 3)\ast n$, we can map the first $b = 1$ layer and the last $a = 3$ layers of $Q_m$ into $Q_{[S,[N]]}$. By Lemma 1, $Q_{[S,[N]]}$ contains either a blue $Q_n$ or a red $Q_m$, so $Q_N$ contains either a blue $Q_n$ or a red $Q_m$. Subcase 4 is similar; we consider $Q_{[S, [N]]}$.





In subcase 2, we consider $Q_{[\emptyset,T]}$, a poset of dimension at least $N - 2$. Both $S$ and $T$ are blue. We consider red sets of maximum and minimum cardinality in $Q_{[\emptyset,T]}$, and apply the same argument we used to get \eqref{eq:blueblue}.
Since $N - 2 > (m - 2 + \frac{9m - 9 }{(2m - 3)(m + 1)})n + m$, $Q_{[\emptyset, T]}$ contains either a blue $Q_n$ or a red $Q_3$, so $Q_N$ contains either a blue $Q_n$ or a red $Q_3$.





In subcase 1, the top 3 layers of $Q_N$ are red. Let $S$ be a element such that $|S| \le |T|$ for all red elements $T$. We consider $Q_{[S,[N]]}$. Since $N > (m - 2 + \frac{9m - 9 }{(2m - 3)(m + 1)})n + m + 2 \ge (m - 2)n +  \frac{3n}{m(2m - 3)} + 5$ for sufficiently large $m$ and $n$, we have

\[
N > (m - 2)n + \frac{3n}{m(2m - 3)} + 5 \ge (m - 2)n + \lceil 1 + \frac{3n}{m (2m - 3)} \rceil + 3
\]
\[
\ge (m - 3)(n + 1) + \ell + 3 = (m - 3) + (m - 3)n + \ell + 3
\]
\[
= \ell + m + (m + 1 - 1 - 3) \ast n.
\]

We map the first $b = 1$ layer and the last $a = 3$ layers of $Q_3$ into $Q_{[S,[N]]}$. By Lemma 1, $Q_{[S,[N]]}$ contains either a blue $Q_n$ or a red $Q_3$, so $Q_N$ contains either a blue $Q_n$ or a red $Q_3$.

\end{proof}

\begin{proof}[\textbf{Proof of Theorem \ref{thm:2 and 3}}]
Consider a coloring $c$ of $Q_4$ defined by

\[ c(S) = \left\{
\begin{array}{ll}
      \text{blue} & \text{if } |S| \text{ is even } \\
      \text{red} & \text{if } |S| \text{ is odd } \\
\end{array}
\right. \]

for all sets $S$ in $Q_4$. This coloring of $Q_4$ contains no red copy of $Q_2$ and no blue copy of $Q_3$. Thus, $R(Q_2, Q_3) > 4$. Now we need only show $R(Q_2, Q_3) \le 5$.

Consider a red-blue coloring of $Q_5$ containing no red $Q_2$ and no blue $Q_3$. We consider the following cases.
\bigskip

\noindent\textbf{Case 1.} \textit{ Both $\emptyset$ and $[5]$ are colored red.}
\bigskip

Let $u,v$ be two red elements in $Q_5$. If $u$ and $v$ are incomparable, $\{ \emptyset, u, v, [5] \}$ form a red $Q_2$. So every red elements in $Q_5$ lies on the same maximal chain. With the exception of this maximal chain, the rest of $Q_5$ is blue, and we can find a blue $Q_3$, a contradiction.

\bigskip

\noindent\textbf{Case 2.} \textit{ One of $\emptyset$ and $[5]$ is colored red, and the other is blue.}
\bigskip

Without loss of generality, suppose $\emptyset$ is red and $[5]$ is blue. Suppose there exists a red set $T$ with $|T| = 4$. Without loss of generality, let $T$ be $\{1, 2, 3, 4 \}$. Consider $Q_{[\emptyset, T]}$, and let $U,V$ be two red elements in $Q_{[\emptyset, T]}$. If $U$ and $V$ are incomparable, $\{\emptyset, U, V, T \}$ form a red $Q_2$. So every red element in $Q_{[\emptyset, T]}$ lies on the same maximal chain. Without loss of generality, suppose this maximal chain is $\{ \emptyset, \{ 1 \}, \{ 1, 2 \}, \{ 1, 2 , 3 \}, \{ 1, 2, 3, 4 \} \}$. Then the sets $\{ 4 \}, \{1,4\}, \{2,4\}, \{3,4\}, \{1,2,4\}, \{1,3,4\}, \{2,3,4\}$ all must be blue. These sets, along with $[5]$ form a blue $Q_3$. Thus, every set in level 4 of $Q_5$ must be blue.

Suppose there exists two red sets $S_1$ and $S_2$ with $|S_1| = |S_2| = 1$. Then $S_1 \cup S_2$ must be blue. Moreover, every set in $Q_{[S_1 \cup S_2, [5]]}$ must be blue. Then $Q_{[S_1 \cup S_2, [5]]}$ is a blue copy of $Q_3$, a contradiction. Thus, $Q_5$ has at most one red level 1 set.

Without loss of generality, suppose $\{1\}$ is the only red level 1 set in $Q_5$. Note that $\bar{2}, \bar{3},$ and $\bar{4}$ are all blue. Consider $Q_{[\{5\}, \bar{2}\cap\bar{3}]}$. If $\bar{2} \cap \bar{3} = \{ 1, 4, 5 \}$ and $\{4 , 5 \}$ are both red, then $\{ \emptyset, \{1\}, \{4,5 \}, \{1,4,5\} \}$ is a red copy of $Q_2$. Thus, at least one of $\{4,5\}$ and $\{1,4,5\}$ is blue. Similarly, when we consider $Q_{[\{5\}, \bar{2}\cap\bar{4}]}$ and $Q_{[\{5\}, \bar{3}\cap\bar{4}]}$, we conclude that at least one of $\{3,5\}$ and $\{1,3,5\}$ is blue and at least one of $\{2,5\}$ and $\{1,2,5\}$ is blue. These blue sets, along with $\{5\}$, $\bar{2}$, $\bar{3}$, $\bar{4}$, and $[5]$ form a blue copy of $Q_3$. Thus, $Q_5$ has no red level 1 set.

Now, note that $\{1,2\}$, $\{1,3\}$, and $\{1,4\}$ cannot all be blue. Otherwise, 

$\{ \{1\}, \{1,2\}, \{1,3\}, \{1,4\}, \bar{2}, \bar{3}, \bar{4}, [5] \}$ is a blue copy of $Q_3$. Suppose, without loss of generality, that $\{1,2\}$ is red. Consider $Q_{[\{1\}, \{1,2,3\}]}$. If $\{2,3\}$ and $\{1,2,3\}$ are both red, then $\{ \emptyset, \{1,2\}, \{2,3\}, \{1,2,3\}\}$ is a red copy of $Q_2$. Thus, at least one of $\{2,3\}$ and $\{1,2,3\}$ is blue. Similarly, when we consider $Q_{[\{1\}, \{1,2,4\}]}$ and $Q_{[\{1\}, \{1,2,5\}]}$, we conclude that at least one of $\{2,4\}$ and $\{1,2,4\}$ is blue and at least one of $\{2,5\}$ and $\{1,2,5\}$ is blue. These blue sets, along with $\{1\}, \bar{3}, \bar{4}, \bar{5},$ and $[5]$ form a blue copy of $Q_3$, a contradiction.



\bigskip

\noindent\textbf{Case 3.} \textit{ Both $\emptyset$ and $[5]$ are colored blue.}
\bigskip

Suppose $Q_5$ has at most 2 red level 1 sets. In other words, $Q_5$ has at least 3 blue level 1 sets. Without loss of generality, suppose $\{1\}$, $\{2\}$, and $\{3\}$ are all blue. Consider $Q_{[\{1,2\}, \bar{3}]}$. If every set in $Q_{[\{1,2\}, \bar{3}]}$ is red, $Q_{[\{1,2\}, \bar{3}]}$ is a red copy of $Q_2$. Thus, there is at least one blue set in $Q_{[\{1,2\}, \bar{3}]}$. Similarly, there is at least one blue set in $Q_{[\{1,3\}, \bar{2}]}$ and at least one blue set in $Q_{[\{2,3\}, \bar{1}]}$. These sets, along with $\emptyset, \{1\}, \{2\}, \{3\},$ and $[5]$, form a blue copy of $Q_3$. Thus, $Q_5$ has at least 3 red level 1 sets. By a similar argument, $Q_5$ also has at least 3 red level 4 sets.

Let $S_1, S_2, S_3$ be 3 red level 1 sets, and let $T_1, T_2, T_3$ be 3 level 4 sets. We consider the following subcases.

\bigskip

\noindent\textbf{Subcase 3.1} \textit{ At least one of $S_1, S_2,$ and $S_3$ is a subset of $T_1, T_2,$ and $T_3$}.
\bigskip

Without loss of generality, let $S_1 = \{1\}$ be red and a subset of $T_1 = \bar{3} = \{1,2,4,5\}$, $T_2 = \bar{4} = \{1,2,3,5\}$, and $T_3 = \bar{5} = \{1,2,3,4\}$, all of which are red. Note that no two of $\{1,2,3 \}$, $\{1,2,4\}$, and $\{1,2,5\}$ can be red without creating a red copy of $Q_2$. 
Also, no two of $\{1,2\}$, $\{1,3\}$, $\{1,4\}$, and $\{1,5\}$ can be red without creating a red copy of $Q_2$.

Suppose $\{1,2\}$ is red, which means $\{1,3\}$, $\{1,4\}$, and $\{1,5\}$ must all be blue, and $\{1,4,5\}$, $\{1,3,5\}$, and $\{1,3,4\}$ must all be blue. These 6 sets, along with $\emptyset$ and $[5]$, form a blue copy of $Q_3$, a contradiction.

Suppose exactly one of $\{1,2,3 \}$, $\{1,2,4\}$, and $\{1,2,5\}$ is red. Without loss of generality, suppose $\{1,2,3\}$ is red. Neither $\{1,4\}$ nor $\{1,5\}$ can be red without creating a red copy of $Q_2$ with $\{1\}$, $\{1,2,3\}$, and $\{1,2,3,4\}$. Suppose $\{1,3\}$ is red, which means $\{1,2\}$, $\{1,4\}$, and $\{1,5\}$ must all be blue. Then $\{1,4,5\}$ must be red. If $\{1,3,4,5\}$ is red, it forms a red copy of $Q_2$ with $\{1\}$, $\{1,3\}$, and $\{1,4,5\}$. If $\{1,3,4,5\}$ is blue, it forms a blue copy of $Q_3$ with $\emptyset$, $\{1,2\}$, $\{1,4\}$, $\{1,5\}$, $\{1,2,4\}$, $\{1,2,5\}$, and $[5]$. Thus, $Q_5$ contains a red copy of $Q_2$ or a blue copy of $Q_3$, a contradiction.

Suppose $\{1,2,3\}$ is red and none of $\{1,3\}$, $\{1,4\}$, and $\{1,5\}$ are red. Then $\{1,4,5\}$ must be red, and $\{1,3,4\}$ and $\{1,3,5\}$ must be blue. Then $\{2,4\}$, $\{2,5\}$, $\{3,4\}$, and $\{3,5\}$ must be blue, and $\{2,4,5\}$ must be red. Then $\{4\}$, $\{5\}$, and $\{4,5\}$ must be blue. Then $\{4\}$, $\{5\}$, $\{1,2\}$, $\{4,5\}$, $\{1,2,4\}$, and $\{1,2,5\}$, along with $\emptyset$ and $[5]$, form a blue copy of $Q_3$, a contradiction.

Now suppose none of $\{1,2,3 \}$, $\{1,2,4\}$, or $\{1,2,5\}$ are red. Again, no two of $\{1,2\}$, $\{1,3\}$, $\{1,4\}$ and $\{1,5\}$ are red. 
Suppose one of $\{1,3\}$, $\{1,4\}$, and $\{1,5\}$ is red. Without loss of generality, suppose $\{1,3\}$ is red. Then $\{1,4,5\}$ is must be red, and $\{1,3,4,5\}$ must be blue. Then $\{1,2\}$, $\{1,4\}$, $\{1,5\}$, $\{1,2,3\}$, $\{1,2,5\}$, and $\{1,3,4,5\}$, along with $\emptyset$ and $[5]$, form a blue copy of $Q_3$, a contradiction.

Suppose none of $\{1,2\}$, $\{1,3\}$, $\{1,4\}$, or $\{1,5\}$ are red. Then $\{1,4,5\}$, $\{1,3,4\}$, and $\{1,3,5\}$ must all be red, and $\{1,3,4,5\}$ must be blue. Then $\{1,2\}$, $\{1,3\}$, $\{1,4\}$, $\{1,2,3\}$, $\{1,2,4\}$, and $\{1,3,4,5\}$, along with $\emptyset$ and $[5]$, form a blue copy of $Q_3$, a contradiction.

In any case where at least one of $S_1,$ $S_2,$ and $S_3$ is a subset of $T_1$, $T_2$, and $T_3$, $Q_5$ contains a red copy of $Q_2$ or a blue copy of $Q_3$.

\bigskip

\noindent\textbf{Subcase 3.2} \textit{ None of $S_1, S_2,$ and $S_3$ is a subset of $T_1, T_2,$ and $T_3$}.
\bigskip

Without loss of generality, let $S_1 = \{1\}$, $S_2 = \{2\}$, $S_3 = \{3\}$, $T_1 = \bar{1} = \{2,3,4,5\}$, $T_2 = \bar{2} = \{1,3,4,5\}$, and $T_3 = \bar{3} = \{1,2,4,5\}$ all be red. Certainly, if every level 2 set and every level 3 set is blue, or if one or both of $\{4,5\}$ and $\{1,2,3\}$ are the only red sets, then $Q_5$ contains a blue copy of $Q_3$.

Suppose one of $\{1,2\}$, $\{1,3\}$ and $\{2,3\}$ is red. Without loss of generality, suppose $\{1,2\}$ is red. Then $\{1,4\}$, $\{1,5\}$, $\{2,4\}$, $\{2,5\}$, $\{1,4,5\}$, and $\{2,4,5\}$ must all be blue. Suppose either $\{1,2,3,4\}$ or $\{1,2,3,5\}$ is red. Without loss of generality, suppose $\{1,2,3,4\}$ is red. Then $\{1,3\}$, $\{2,3\}$, $\{1,3,4\}$, and $\{2,3,4\}$ must all be blue, and $\{1,2,3,5\}$ must be red. Then $\{1,3,5\}$ and $\{2,3,5\}$ must be blue. The sets $\{1,4\}$, $\{1,5\}$, $\{1,3\}$, $\{1,4,5\}$, $\{1,3,4\}$, and $\{1,3,5\}$, along with $\emptyset$ and $[5]$, form a blue copy of $Q_3$, a contradiction.

Now suppose $\{1,2\}$ is red and $\{1,2,3,4\}$ and $\{1,2,3,5\}$ are both blue. Then $\{1,3\}$ must be red, and $\{1,2,3\}$ must be blue. Then $\{1,4\}$, $\{1,5\}$, $\{1,2,3\}$, $\{1,4,5\}$, $\{1,2,3,4\}$, and $\{1,2,3,5\}$, along with $\emptyset$ and $[5]$, form a blue copy of $Q_3$, a contradiction. The argument is similar if any one of $\{1,4,5\}$, $\{2,4,5\}$, and $\{3,4,5\}$ is red.

Suppose any level 2 set other than $\{1,2\}$, $\{1,3\}$, $\{2,3\}$, or $\{4,5\}$ is red. Without loss of generality, suppose $\{1,4\}$ is red. Then $\{1,2\}$, $\{1,3\}$, $\{1,5\}$, $\{1,2,5\}$, and $\{1,3,5\}$ are all blue. Then $\{1,2,3\}$ must be red, and $\{1,2,3,4\}$ must be blue. Then $\{1,2\}$, $\{1,3\}$, $\{1,5\}$, $\{1,2,5\}$, $\{1,3,5\}$, and $\{1,2,3,4\}$, along with $\emptyset$ and $[5]$, form a blue copy of $Q_3$, a contradiction. The argument is similar if any level 3 set other than $\{1,4,5\}$, $\{2,4,5\}$, $\{3,4,5\}$, or $\{1,2,3\}$ is red.

In any case where none of $S_1$, $S_2$, and $S_3$ is a subset of $T_1$, $T_2$, and $T_3$, $Q_5$ contains a red copy of $Q_2$ or a blue copy of $Q_3$.

\end{proof}

\section{Concluding Remarks}\label{remarks} 
There remains a significant gap between our upper bounds and the best known lower bounds given by Axenovich and Walzer. We believe the true values of $R(Q_m, Q_n)$ for sufficiently large $m$ and $n$ are significantly less than our upper bounds. Assuming, without loss of generality, that $n \ge m$, we make the following conjecture for sufficiently large $m$ and $n$.

\begin{conjecture}
$R(Q_m, Q_n) = o(n^2)$.
\end{conjecture}

\end{document}